\documentclass[12pt]{article}

\usepackage{a4wide}

\usepackage[intlimits]{amsmath}
\usepackage{amssymb,amsthm}

\usepackage[cp1251]{inputenc}
\usepackage[english]{babel}

\theoremstyle{definition}

\theoremstyle{remark}

\numberwithin{equation}{section}

\begin{document}\Large

\title{DISTRIBUTION FUNCTION OF MARKOVIAN RANDOM EVOLUTION IN $R^n$}

\author{I.V.Samoilenko\\
Institute of Mathematics,\\ Ukrainian National Academy of
Sciences,\\ 3 Tereshchenkivs'ka, Kyiv, 01601, Ukraine\\
isamoil@imath.kiev.ua}

\maketitle

\abstract{Obvious view of distribution function of Markovian
random evolution is found in terms of Bessel functions of $n+1$-th
order. 
}

\section{\bf  Introduction.} Markovian random evolutions in $R^n$ were
studied in the work [1], where the connection between the equation
for the functional of the evolution and Bessel equation of high
order was found.

Bessel equations of high order and their solutions - Bessel
functions of high order - were studied by Turbin and Plotkin in
[2].

In the works [3,4] Orsingher and Sommella managed to combine these
two results and to receive the obvious view for the distribution
functions of random evolutions in $R^2$ and $R^3$. In [4] they
also made a conjuncture about the view of the distribution
function of the random evolution in $R^n$.

In this paper we generalize the results of Orsingher and Sommella
and prove their conjuncture.

Let us start with the definition of Markovian random evolution.To
do this we have to define a regular $n+1$-hedron inscribed into a
unit sphere in $R^n$.

{\bf Definition 1.1:} We call figure in $R^n$ a regular
$n+1$-hedron inscribed into a unit sphere if the following
conditions are true:

1. The figure has $n+1$ vertices $\tau_0,\ldots,\tau_n$ situated
at the unit sphere with the center in $(0,\ldots,0)$.

2. The center of masses of the figure is situated in
$(0,\ldots,0)$.

3. The $i$-th vertice is situated into $i+1$-dimensional subspace,
determined by the first $i+1$ coordinates of $R^n$.

Having this definition, we may easily find the coordinates of the
vertices. The components of ${\tau}_i, i=\overline{0,n}$ are equal
to
$$\tau_i^j=\left\{\begin{array}{c}
  -\frac{1}{n}\sqrt{\frac{n(n+1)}{(n-j+1)(n-j+2)}}, j<i+1 \\
  \sqrt{\frac{(n+1)(n-i)}{n(n-i+1)}}, j=i+1 \\
  0, j>i+1 \\
\end{array}\right., j=\overline{1,n}. \eqno(1)$$

Really, the first condition is true, because the distance from
every vertice to the zero point is equal to 1.

The second condition is true, because the sum of first components
of all coordinates is equal to 0, the same for second components,
etc.

At last, for the $i$-th vertice $i+2,i+3,\ldots$ components of
coordinates are equal to 0. So the third condition is also true.

We are now ready to define Markovian random evolution in $R^n$
from [1].

{\bf Definition 1.2:} We call a random process $\overline{S}(t)$
Markovian random evolution in $R^n$ if:

1. The process begins its motion at the point
$\overline{x}=(x_1,\ldots,x_n).$

2. The initial direction of the motion is $\overline{\tau}_i,
i=\overline{0,n},$ where the components of coordinates of
$\overline{\tau}_i$ are presented in (1).

3. The time of the motion at some direction is distributed like
$e^{-\lambda t}.$

4. The $k$-th direction is followed by the $k+1$-th direction (the
$n$-th one is followed by 0).

5. The velocity of the particle's motion is equal to $v$.

The system of backward Kolmogorov equations and corresponding
$n+1$-dimentional equation for the functionals (here
$i=\overline{0,n}$ is the start direction)
$$u_i(\overline{x},t)=E_i(f(\overline{S}(t)) \eqno(2)$$
were also studied in [1].

Using a well-known result from the theory of Markovian processes
we may receive the system of differential equations for the
distribution functions $$f_j(x_1,\ldots,x_n,t)dx_1\ldots
dx_n=P\{S_1(t)\in dx_1,\ldots,S_n(t)\in dx_n,D(t)=j\}, \eqno(3)$$
where D(t) is the direction at time $t$.

The system of equations for the distribution functions (3) is
adjoint to the system for the functionals (2) (see [1]):
$$ \begin{cases}
    \frac{\partial}{\partial t}f_0(\overline{x},t)=-\lambda f_0(\overline{x},t)-v
    (\overline{\tau}_0,\nabla)f_0(\overline{x},t)+\lambda f_n(\overline{x},t)\\
    \frac{\partial}{\partial t}f_1(\overline{x},t)=-\lambda f_1(\overline{x},t)-v
    (\overline{\tau}_1,\nabla)f_1(\overline{x},t)+\lambda f_0(\overline{x},t)\\
    \ldots \\
    \frac{\partial}{\partial t}f_n(\overline{x},t)=-\lambda f_n(\overline{x},t)-v
    (\overline{\tau}_n,\nabla)f_n(\overline{x},t)+\lambda f_{n-1}(\overline{x},t)
  \end{cases} \eqno(4)$$

  The corresponding $n+1$-dimentional equation is: $$\prod_{i=0}^n\left(\frac{\partial}{\partial
t}+\lambda+v(\overline{\tau}_i,\nabla)\right)f(\overline{x},t)=\lambda^{n+1}f(\overline{x},t),\eqno(5)$$
where $f$ is any of functions $f_0,\ldots,f_n$.

To describe the support of the distribution of random evolution we
use the notion of regular $n+1$-hedron from {\bf Definition 1.1}.

One of the most important characteristics of random evolution that
differs it from Brownian motion is that evolution with probability
1 moves into a regular $n+1$-hedron $$T_{vt}\equiv
\left\{x_1,\ldots,x_n: -\frac{vt}{n}<x_1<vt,
\frac{1}{n-k+1}\left[\sqrt{\frac{(n-1)\ldots (n-k+1)}{(n+1)\ldots
(n-k+3)}}x_1+\right.\right.$$
$$\left.+\ldots+\sqrt{\frac{(n-k+1)}{(n-k+3)}}x_{k-1}-
vt\sqrt{\frac{(n-1)\ldots (n-k+1)}{(n+1)\ldots (n-k+3)}}
\right]<x_k<$$
$$<-\left[\sqrt{\frac{(n-1)\ldots (n-k+1)}{(n+1)\ldots
(n-k+3)}}x_1+\ldots+\sqrt{\frac{(n-k+1)}{(n-k+3)}}x_{k-1}-\right.$$
$$\left.- vt\sqrt{\frac{(n-1)\ldots (n-k+1)}{(n+1)\ldots (n-k+3)}},
k=\overline{2,n} \right\}.\eqno(6)$$

The probability that evolution is on the edge $\partial T_{vt}$ of
$n+1$-hedron is equal to: $$P\{\overline{S}(t)\in \partial
T_{vt}\}=e^{-\lambda t}+\lambda te^{-\lambda t}+\ldots
+\frac{(\lambda t)^{n-1}}{(n-1)!}e^{-\lambda t}.$$ Here
$e^{-\lambda t}$ is the probability of being at time $t$ on the
vertices of $T_{vt}$, $\frac{(\lambda t)^{k-1}}{(k-1)!}e^{-\lambda
t}, k=\overline{1,n-1}$ is the probability of being on
$k$-dimensional edge of $T_{vt}$.

So, the continuous part of the distribution we are to find has the
property:
$$\int\ldots\int_{T_{vt}}f(x_1,\ldots,x_n;t)dx_1\ldots dx_n=1-e^{-\lambda t}-\lambda te^{-\lambda t}-\ldots
-\frac{(\lambda t)^{n-1}}{(n-1)!}e^{-\lambda t}. \eqno(7)$$

In other words, we have to find a non-negative continuous function
that satisfies system (4) and condition (7).

\section{\bf  Derivation of the equations.}

Making exponential transformation $f_j=e^{-\lambda t}g_j,
j=\overline{0,n}$ we receive from (4):$$\begin{cases}
    \frac{\partial}{\partial t}g_0(\overline{x},t)=-v
    (\overline{\tau}_0,\nabla)g_0(\overline{x},t)+\lambda g_n(\overline{x},t)\\
    \frac{\partial}{\partial t}g_1(\overline{x},t)=-v
    (\overline{\tau}_1,\nabla)g_1(\overline{x},t)+\lambda g_0(\overline{x},t)\\
    \ldots \\
    \frac{\partial}{\partial t}g_n(\overline{x},t)=-v
    (\overline{\tau}_n,\nabla)g_n(\overline{x},t)+\lambda g_{n-1}(\overline{x},t)
  \end{cases}\eqno(8)$$

  {\bf Theorem 2.1} The functions $g_j, j=\overline{0,n}$ are
  solutions of the equation $$\frac{\partial g}{\partial y_1\ldots\partial y_{n+1}}=
  \left(\frac{\lambda}{v}\right)^{n+1}\left(\frac{\sqrt{n}}{\sqrt{n+1}}\right)^{n+1}\left(\frac{1}{\sqrt[2n+2]{2n+2}}\right)^{n+1}g,\eqno(9)$$
  where $$y_1=\frac{vt}{n}+x_1,$$ $$y_k=\frac{1}{n-k+1}\left[\sqrt{\frac{(n-1)\ldots (n-k+1)}{(n+1)\ldots
  (n-k+3)}}vt-\sqrt{\frac{(n-1)\ldots (n-k+1)}{(n+1)\ldots
  (n-k+3)}}x_1-\right.$$ $$\left.-\ldots-\sqrt{\frac{n-k+1}{n-k+3}}x_{k-1}\right]+x_k,
  k=\overline{2,n},$$ $$y_{n+1}=\sqrt{\frac{(n-1)\ldots 1}{(n+1)\ldots 3}}vt-\sqrt{\frac{(n-1)\ldots 1}{(n+1)\ldots 3}}x_1
  -\ldots-\sqrt{\frac{1}{3}}x_{n-1}-x_n.\eqno(10)$$

  {\it Proof.} By applying the transformation (10) to the system
  (8) we receive: $$\begin{cases}
    \frac{n+1}{n}v\frac{\partial g_0}{\partial y_1}=\lambda g_n\\
    \frac{\sqrt{(n+1)(n-k+1)}}{\sqrt{n(n-k)}}v \frac{\partial g_k}{\partial y_{k+1}}=\lambda g_{k-1}, k=\overline{1,n-1}\\
    \frac{\sqrt{2(n+1)}}{\sqrt{n}}v\frac{\partial g_n}{\partial y_{n+1}}=\lambda
    g_{n-1}.
  \end{cases}$$ By differentiation we easily have (9).

  Theorem is proved.

  {\bf Theorem 2.2} The transformation $z=\sqrt[n+1]{y_1\ldots
  y_{n+1}}$converts the equation (9) into $n+1$-dimensional Bessel
  equation: $$\left(z\frac{\partial }{\partial z}\right)^{n+1}g=\left(\frac{\lambda}{v}\right)^{n+1}\left(\sqrt{n(n+1)}\right)^{n+1}
  \left(\frac{1}{\sqrt[2n+2]{2n+2}}\right)^{n+1}z^{n+1}g, \eqno(11)$$ where
  $g$ is any of the functions $g_j, j=\overline{0,n}.$

  {\it Proof.} We don't show the proof since plain calculations are
  involved.

  The solutions of equation (11) were found in [2]. We use
  here one of the functions that is solution of (11): $$I_{0,n}\left(\frac{\lambda}{v}
  \frac{\sqrt{n(n+1)}}{\sqrt[2n+2]{2n+2}}z\right)=\sum_{k=0}^{\infty}\left(\frac{\lambda}{v}
  \frac{\sqrt{n(n+1)}}{\sqrt[2n+2]{2n+2}}\frac{1}{n+1}z\right)^{(n+1)k}\frac{1}{(k!)^{n+1}}.\eqno(12)$$

  \underline{Remark 2.1} The proofs that this function is solution
  of (11) in case of $n=2,3$ may also be found in [3,4].

  Making the backward change of variables in (12) and exponential transformation we receive
  the solution of system (4).

\section{\bf Probability distribution.}

  Now we are ready to formulate the main result of the article:

  {\bf Theorem 3.1} The absolutely continuous component of the
  distribution of Markovian random evolution in $R^n$ is: $$f(x_1,\ldots,x_n;t)
  =\frac{(\sqrt{n})^ne^{-\lambda t}}{(\sqrt{n+1})^{n+1}v^n}\left[\lambda^n+\lambda^{n-1}
  \frac{\partial}{\partial t}+\ldots+\right.$$ $$\left.+\frac{\partial^n}{\partial t^n}\right]I_{0,n}\left(\frac{\lambda}{v}
  \frac{\sqrt{n(n+1)}}{\sqrt[2n+2]{2n+2}}\sqrt[n+1]{y_1\ldots
  y_{n+1}}\right),\eqno(13)$$ here and later $y_i, i=\overline{1,n+1}$ are
  the linear combinations of $t,x_i, i=\overline{1,n}$
  pointed in (10).

  {\it Proof.} The fact that (13) is a non-negative function may be
  easily seen from the view of the series (12).

  The function (13) satisfies the equation (5).Really, $I_{0,n}\left(\frac{\lambda}{v}
  \frac{\sqrt{n(n+1)}}{\sqrt[2n+2]{2n+2}}\times \right.$ $\left.\times\sqrt[n+1]{y_1\ldots
  y_{n+1}}\right),$ where $y_i, i=\overline{1,n+1}$ are pointed
  in (10) satisfies equation (9). Any linear combination of
  derivatives of the last function with
  the coefficient $\frac{(\sqrt{n})^n}{(\sqrt{n+1})^{n+1}v^n}$ also satisfies (9).
  So, if we make backward exponential transformation, we receive a
  solution of (5).

  The last problem is to prove the condition (7) for the function $f(x_1,\ldots,x_n;t)$.

  We first note that $$\int\ldots\int_{T_{vt}}I_{0,n}\left(\frac{\lambda}{v}
  \frac{\sqrt{n(n+1)}}{\sqrt[2n+2]{2n+2}}\sqrt[n+1]{y_1\ldots
  y_{n+1}}\right)dx_1\ldots dx_n=$$ $$=\sum_{k=0}^{\infty}\left(\frac{\lambda}{(n+1)v}\right)^{(n+1)k}
  \frac{(\sqrt{n(n+1)})^{(n+1)k}}{(\sqrt[2n+2]{2n+2})^{(n+1)k}}\frac{1}{(k!)^{n+1}}
  \int_{-\frac{vt}{n}}^{vt}\frac{1}{n^k}[vt+nx_1]^kdx_1\ldots$$ $$\int_{\frac{1}{2}\left[\sqrt{
  \frac{(n-1)\ldots 2}{(n+1)\ldots 4}}x_1+\ldots+\sqrt{\frac{2}{4}}x_{n-2}-\sqrt{
  \frac{(n-1)\ldots 2}{(n+1)\ldots 4}}vt\right]}^{-\left[\sqrt{
  \frac{(n-1)\ldots 2}{(n+1)\ldots 4}}x_1+\ldots+\sqrt{\frac{2}{4}}x_{n-2}-\sqrt{
  \frac{(n-1)\ldots 2}{(n+1)\ldots 4}}vt\right]}\frac{1}{2^k}\left[-\sqrt{
  \frac{(n-1)\ldots 2}{(n+1)\ldots 4}}x_1-\ldots-\right.$$ $$\left.-\sqrt{\frac{2}{4}}x_{n-2}+\sqrt{
  \frac{(n-1)\ldots 2}{(n+1)\ldots
  4}}vt+2x_{n-1}\right]^kdx_{n-1}\times$$ $$\times
  \int_{\left[\sqrt{
  \frac{(n-1)\ldots 1}{(n+1)\ldots 3}}x_1+\ldots+\sqrt{\frac{1}{3}}x_{n-1}-\sqrt{
  \frac{(n-1)\ldots 1}{(n+1)\ldots 3}}vt\right]}^{-\left[\sqrt{
  \frac{(n-1)\ldots 1}{(n+1)\ldots 3}}x_1+\ldots+\sqrt{\frac{1}{3}}x_{n-1}-\sqrt{
  \frac{(n-1)\ldots 1}{(n+1)\ldots 3}}vt\right]}\left[\left(-\sqrt{
  \frac{(n-1)\ldots 1}{(n+1)\ldots 3}}x_1-\right.\right.$$ $$\left.\left.-\ldots-\sqrt{\frac{1}{3}}x_{n-1}+\sqrt{
  \frac{(n-1)\ldots 1}{(n+1)\ldots 3}}vt\right)^2-x_{n}^2\right]^kdx_{n}.\eqno(14)$$

  The inner, first, integral could be found using the change of variables
  $x_n=\left(-\sqrt{
  \frac{(n-1)\ldots 1}{(n+1)\ldots 3}}x_1-\ldots-\sqrt{\frac{1}{3}}x_{n-1}+\sqrt{
  \frac{(n-1)\ldots 1}{(n+1)\ldots 3}}vt\right)z$. We have:$$\int_{\left[\sqrt{
  \frac{(n-1)\ldots 1}{(n+1)\ldots 3}}x_1+\ldots+\sqrt{\frac{1}{3}}x_{n-1}-\sqrt{
  \frac{(n-1)\ldots 1}{(n+1)\ldots 3}}vt\right]}^{-\left[\sqrt{
  \frac{(n-1)\ldots 1}{(n+1)\ldots 3}}x_1+\ldots+\sqrt{\frac{1}{3}}x_{n-1}-\sqrt{
  \frac{(n-1)\ldots 1}{(n+1)\ldots 3}}vt\right]}\left[\left(-\sqrt{
  \frac{(n-1)\ldots 1}{(n+1)\ldots 3}}x_1-\ldots-\sqrt{\frac{1}{3}}x_{n-1}+\right.\right.$$ $$\left.\left.+\sqrt{
  \frac{(n-1)\ldots 1}{(n+1)\ldots 3}}vt\right)^2-x_{n}^2\right]dx_{n}=\left(-\sqrt{
  \frac{(n-1)\ldots 1}{(n+1)\ldots 3}}x_1-\ldots-\sqrt{\frac{1}{3}}x_{n-1}+\right.$$ $$\left.+\sqrt{
  \frac{(n-1)\ldots 1}{(n+1)\ldots 3}}vt\right)^{2k+1}\int_{-1}^1(1-z^2)^kdz=\left(-\sqrt{
  \frac{(n-1)\ldots 1}{(n+1)\ldots 3}}x_1-\ldots-\right.$$ $$\left.-\sqrt{\frac{1}{3}}x_{n-1}+\sqrt{\frac{1}{3}}x_{n-1}+\sqrt{
  \frac{(n-1)\ldots 1}{(n+1)\ldots 3}}vt\right)^{2k+1}\frac{2^{2k+1}(k!)^2}{(2k+1)!}.$$

  In the following, second, integral $$\int_{\frac{1}{2}\left[\sqrt{
  \frac{(n-1)\ldots 2}{(n+1)\ldots 4}}x_1+\ldots+\sqrt{\frac{2}{4}}x_{n-2}-\sqrt{
  \frac{(n-1)\ldots 2}{(n+1)\ldots 4}}vt\right]}^{-\left[\sqrt{
  \frac{(n-1)\ldots 2}{(n+1)\ldots 4}}x_1+\ldots+\sqrt{\frac{2}{4}}x_{n-2}-\sqrt{
  \frac{(n-1)\ldots 2}{(n+1)\ldots 4}}vt\right]}\frac{1}{2^k}\left[-\sqrt{
  \frac{(n-1)\ldots 2}{(n+1)\ldots 4}}x_1-\ldots-\sqrt{\frac{2}{4}}x_{n-2}+\right.$$ $$\left.+\sqrt{
  \frac{(n-1)\ldots 2}{(n+1)\ldots 4}}vt+2x_{n-1}\right]^k\left[-\sqrt{
  \frac{(n-1)\ldots 1}{(n+1)\ldots 3}}x_1-\ldots-\sqrt{\frac{1}{3}}x_{n-1}+\right.$$ $$\left.+\sqrt{
  \frac{(n-1)\ldots 1}{(n+1)\ldots 3}}vt\right]^{2k+1}dx_{n-1}=\int_{\frac{1}{2}\left[\sqrt{
  \frac{(n-1)\ldots 2}{(n+1)\ldots 4}}x_1+\ldots+\sqrt{\frac{2}{4}}x_{n-2}-\sqrt{
  \frac{(n-1)\ldots 2}{(n+1)\ldots 4}}vt\right]}^{-\left[\sqrt{
  \frac{(n-1)\ldots 2}{(n+1)\ldots 4}}x_1+\ldots+\sqrt{\frac{2}{4}}x_{n-2}-\sqrt{
  \frac{(n-1)\ldots 2}{(n+1)\ldots 4}}vt\right]}\frac{1}{2^k\sqrt{3}^{2k+1}}\times$$ $$\times\left[-\sqrt{
  \frac{(n-1)\ldots 2}{(n+1)\ldots 4}}x_1-\ldots-\sqrt{\frac{2}{4}}x_{n-2}+\sqrt{
  \frac{(n-1)\ldots 2}{(n+1)\ldots 4}}vt+\right.$$ $$\left.+2x_{n-1}\right]^k\left[-\sqrt{
  \frac{(n-1)\ldots 2}{(n+1)\ldots 4}}x_1-\ldots-\sqrt{\frac{2}{4}}x_{n-2}+\sqrt{
  \frac{(n-1)\ldots 2}{(n+1)\ldots 4}}vt-\right.$$ $$\left.-x_{n-1}\right]^{2k+1}dx_{n-1}=\frac{1}{2^k\sqrt{3}^{2k+1}}\left(-\sqrt{
  \frac{(n-1)\ldots 2}{(n+1)\ldots 4}}x_1-\ldots-\sqrt{\frac{2}{4}}x_{n-2}+\right.$$ $$\left.+\sqrt{
  \frac{(n-1)\ldots 2}{(n+1)\ldots 4}}vt\right)^{3k+1}\int_{\frac{1}{2}\left[\sqrt{
  \frac{(n-1)\ldots 2}{(n+1)\ldots 4}}x_1+\ldots+\sqrt{\frac{2}{4}}x_{n-2}-\sqrt{
  \frac{(n-1)\ldots 2}{(n+1)\ldots 4}}vt\right]}^{-\left[\sqrt{
  \frac{(n-1)\ldots 2}{(n+1)\ldots 4}}x_1+\ldots+\sqrt{\frac{2}{4}}x_{n-2}-\sqrt{
  \frac{(n-1)\ldots 2}{(n+1)\ldots 4}}vt\right]}(1+$$ $$\left.+2\frac{x_{n-1}}{\left(-\sqrt{
  \frac{(n-1)\ldots 2}{(n+1)\ldots 4}}x_1-\ldots-\sqrt{\frac{2}{4}}x_{n-2}+\sqrt{
  \frac{(n-1)\ldots 2}{(n+1)\ldots 4}}vt\right)}\right)^k\times$$ $$\times\left(1-\frac{x_{n-1}}{\left(-\sqrt{
  \frac{(n-1)\ldots 2}{(n+1)\ldots 4}}x_1-\ldots-\sqrt{\frac{2}{4}}x_{n-2}+\sqrt{
  \frac{(n-1)\ldots 2}{(n+1)\ldots 4}}vt\right)}\right)^{2k+1}dx_{n-1}=$$ $$=\frac{3^{3k+1}}{2^{3k+1}\sqrt{3}^{2k+1}}\left(-\sqrt{
  \frac{(n-1)\ldots 2}{(n+1)\ldots 4}}x_1-\ldots-\sqrt{\frac{2}{4}}x_{n-2}+\sqrt{
  \frac{(n-1)\ldots 2}{(n+1)\ldots 4}}vt\right)^{3k+1}\times$$ $$\times\int_{\frac{1}{2}\left[\sqrt{
  \frac{(n-1)\ldots 2}{(n+1)\ldots 4}}x_1+\ldots+\sqrt{\frac{2}{4}}x_{n-2}-\sqrt{
  \frac{(n-1)\ldots 2}{(n+1)\ldots 4}}vt\right]}^{-\left[\sqrt{
  \frac{(n-1)\ldots 2}{(n+1)\ldots 4}}x_1+\ldots+\sqrt{\frac{2}{4}}x_{n-2}-\sqrt{
  \frac{(n-1)\ldots 2}{(n+1)\ldots 4}}vt\right]}\left(\frac{1}{3}+\frac{2}{3}\times\right.$$ $$\left.\times\frac{x_{n-1}}{\left(-\sqrt{
  \frac{(n-1)\ldots 2}{(n+1)\ldots 4}}x_1-\ldots-\sqrt{\frac{2}{4}}x_{n-2}+\sqrt{
  \frac{(n-1)\ldots 2}{(n+1)\ldots 4}}vt\right)}\right)^k\times$$ $$\times\left(\frac{2}{3}-\frac{2}{3}\frac{x_{n-1}}{\left(-\sqrt{
  \frac{(n-1)\ldots 2}{(n+1)\ldots 4}}x_1-\ldots-\sqrt{\frac{2}{4}}x_{n-2}+\sqrt{
  \frac{(n-1)\ldots 2}{(n+1)\ldots 4}}vt\right)}\right)^{2k+1}dx_{n-1}$$ we make
  the change of variables $z=\frac{1}{3}+\frac{2}{3}\frac{x_{n-1}}{\left(-\sqrt{
  \frac{(n-1)\ldots 2}{(n+1)\ldots 4}}x_1-\ldots-\sqrt{\frac{2}{4}}x_{n-2}+\sqrt{
  \frac{(n-1)\ldots 2}{(n+1)\ldots 4}}vt\right)}$ and so we have $$\frac{3^{3k+2}}{2^{3k+2}\sqrt{3}^{2k+1}}\left(-\sqrt{
  \frac{(n-1)\ldots 2}{(n+1)\ldots 4}}x_1-\ldots-\sqrt{\frac{2}{4}}x_{n-2}+\sqrt{
  \frac{(n-1)\ldots 2}{(n+1)\ldots 4}}vt\right)^{3k+2}\times$$ $$\times\int_0^1z^k(1-z)^{2k+1}dz=\frac{3^{3k+2}}{2^{3k+2}\sqrt{3}^{2k+1}}\left(-\sqrt{
  \frac{(n-1)\ldots 2}{(n+1)\ldots 4}}x_1-\ldots-\sqrt{\frac{2}{4}}x_{n-2}+\right.$$ $$\left.+\sqrt{
  \frac{(n-1)\ldots 2}{(n+1)\ldots 4}}vt\right)^{3k+2}\frac{\Gamma (k+1)\Gamma (2k+2)}{\Gamma (3k+3)}.$$

  The following integrals could be found in the same way. For the
  $m$-th integral we have: $$\left(\sqrt{\frac{m-1}{m+1}}\right)^{m(k+1)-1}\left(\frac{m+1}{m}\right)^{(m+1)(k+1)-1}
  \frac{\Gamma (m(k+1))\Gamma (k+1)}{\Gamma ((m+1)(k+1))}\times$$ $$\times\left(-\sqrt{\frac{(n-1)\ldots m}{(n+1)\ldots
  (m+2)}}x_1-\ldots-\sqrt{\frac{m}{m+2}}x_{n-m}+\right.$$ $$\left.+\sqrt{\frac{(n-1)\ldots m}{(n+1)\ldots
  (m+2)}}vt\right)^{(m+1)(k+1)-1}.$$

  Substituting the integrals found into (14) and performing simple
  calculations we receive: $$\int\ldots\int_{T_{vt}}I_{0,n}\left(\frac{\lambda}{v}
  \frac{\sqrt{n(n+1)}}{\sqrt[2n+2]{2n+2}}\sqrt[n+1]{y_1\ldots
  y_{n+1}}\right)=\left(\frac{v}{\lambda}\right)^n\frac{(\sqrt{n+1})^{n+1}}{(\sqrt{n})^n}\times$$ $$\times\sum_{k=0}^{\infty}
  \frac{(\lambda t)^{(n+1)(k+1)-1}}{((n+1)(k+1)-1)!}.$$

  But to prove (7) we have to find the following integrals:

  $\int\ldots\int_{T_{vt}}\frac{\partial^m}{\partial t^m}
  I_{0,n}\left(\frac{\lambda}{v}
  \frac{\sqrt{n(n+1)}}{\sqrt[2n+2]{2n+2}}\sqrt[n+1]{y_1\ldots
  y_{n+1}}\right)dx_1\ldots dx_n, m=\overline{1,n}.$

  To do this let us find $$\frac{\partial}{\partial t}\int\ldots\int_{T_{vt}}
  I_{0,n}\left(\frac{\lambda}{v}
  \frac{\sqrt{n(n+1)}}{\sqrt[2n+2]{2n+2}}\sqrt[n+1]{y_1\ldots
  y_{n+1}}\right)dx_1\ldots dx_n=\frac{\partial}{\partial t}\int_{-\frac{vt}{n}}^{vt}\ldots$$
  $$\int_{\frac{1}{2}\left[\sqrt{
  \frac{(n-1)\ldots 2}{(n+1)\ldots 4}}x_1+\ldots+\sqrt{\frac{2}{4}}x_{n-2}-\sqrt{
  \frac{(n-1)\ldots 2}{(n+1)\ldots 4}}vt\right]}^{-\left[\sqrt{
  \frac{(n-1)\ldots 2}{(n+1)\ldots 4}}x_1+\ldots+\sqrt{\frac{2}{4}}x_{n-2}-\sqrt{
  \frac{(n-1)\ldots 2}{(n+1)\ldots 4}}vt\right]}
  \int_{\left[\sqrt{
  \frac{(n-1)\ldots 1}{(n+1)\ldots 3}}x_1+\ldots+\sqrt{\frac{1}{3}}x_{n-1}-\sqrt{
  \frac{(n-1)\ldots 1}{(n+1)\ldots 3}}vt\right]}^{-\left[\sqrt{
  \frac{(n-1)\ldots 1}{(n+1)\ldots 3}}x_1+\ldots+\sqrt{\frac{1}{3}}x_{n-1}-\sqrt{
  \frac{(n-1)\ldots 1}{(n+1)\ldots
  3}}vt\right]}I_{0,n}\left(\frac{\lambda}{v}\times\right.$$ $$\left.\times
  \frac{\sqrt{n(n+1)}}{\sqrt[2n+2]{2n+2}}\sqrt[n+1]{y_1\ldots
  y_{n+1}}\right)dx_1\ldots dx_{n-1}dx_{n}.$$

  The main rule for such integrals is: we find the derivative of
  the above limit and multiply it on the integral in which we
  substitute the above limit instead of $x_1$, then we subtract
  the derivative of the down limit multiplied on the integral in which we
  substitute the down limit instead of $x_1$.On so on for every
  integral.

  But in our case when we substitute the above limit into the next
  integral  it becomes equal to 0 - this may be easily seen using
  (6). Otherwise, when we substitute the down limit the argument
  of Bessel function becomes equal to 0 - this may be seen from
  (6) and (10). But it's clear that $I_{0,n}(0)=1$, so the
  integrand will be always equal to 1.

  In other words, we have the following equality: $$\frac{\partial}{\partial t}\int\ldots\int_{T_{vt}}
  I_{0,n}\left(\frac{\lambda}{v}
  \frac{\sqrt{n(n+1)}}{\sqrt[2n+2]{2n+2}}\sqrt[n+1]{y_1\ldots
  y_{n+1}}\right)dx_1\ldots dx_n=$$ $$=\frac{\partial}{\partial t}\int\ldots\int_{T_{vt}}
  1dx_1\ldots dx_n+\int\ldots\int_{T_{vt}}\frac{\partial}{\partial t}
  I_{0,n}\left(\frac{\lambda}{v}
  \frac{\sqrt{n(n+1)}}{\sqrt[2n+2]{2n+2}}\sqrt[n+1]{y_1\ldots
  y_{n+1}}\right)dx_1\ldots dx_n.$$

  But the first integral is the derivative of $Vol T_{vt}.$ We may
  easily find that $Vol T_{vt}=\frac{(\sqrt{n+1})^{n+1}}{(\sqrt{n})^nn!}(vt)^n.$

  So, $$\frac{\partial}{\partial t}\int\ldots\int_{T_{vt}}
  I_{0,n}\left(\frac{\lambda}{v}
  \frac{\sqrt{n(n+1)}}{\sqrt[2n+2]{2n+2}}\sqrt[n+1]{y_1\ldots
  y_{n+1}}\right)dx_1\ldots dx_n=$$ $$=\frac{(\sqrt{n+1})^{n+1}}{(\sqrt{n})^n(n-1)!}v^nt^{n-1}+\int\ldots\int_{T_{vt}}\frac{\partial}{\partial t}
  I_{0,n}\left(\frac{\lambda}{v}
  \frac{\sqrt{n(n+1)}}{\sqrt[2n+2]{2n+2}}\sqrt[n+1]{y_1\ldots
  y_{n+1}}\right)dx_1\ldots dx_n.$$

  For the next derivative we have: $$\frac{\partial^2}{\partial t^2}\int\ldots\int_{T_{vt}}
  I_{0,n}\left(\frac{\lambda}{v}
  \frac{\sqrt{n(n+1)}}{\sqrt[2n+2]{2n+2}}\sqrt[n+1]{y_1\ldots
  y_{n+1}}\right)dx_1\ldots dx_n=$$ $$=\frac{(\sqrt{n+1})^{n+1}}{(\sqrt{n})^n(n-2)!}v^nt^{n-2}+
  \frac{\partial}{\partial t}\int\ldots\int_{T_{vt}}\frac{\partial}{\partial t}
  I_{0,n}\left(\frac{\lambda}{v}
  \frac{\sqrt{n(n+1)}}{\sqrt[2n+2]{2n+2}}\sqrt[n+1]{y_1\ldots
  y_{n+1}}\right)dx_1\ldots dx_n.$$

  We'll find the last integral by the rule we mentioned above, but
  here the above limit makes the following integral equal to 0,
  and when we substitute the down limit into the integrand we have
  $\frac{\partial}{\partial t}
  I_{0,n}(0)=0.$ So, $$\frac{\partial^2}{\partial t^2}\int\ldots\int_{T_{vt}}
  I_{0,n}\left(\frac{\lambda}{v}
  \frac{\sqrt{n(n+1)}}{\sqrt[2n+2]{2n+2}}\sqrt[n+1]{y_1\ldots
  y_{n+1}}\right)dx_1\ldots dx_n=$$ $$=\frac{(\sqrt{n+1})^{n+1}}{(\sqrt{n})^n(n-2)!}v^nt^{n-2}+
  \int\ldots\int_{T_{vt}}\frac{\partial^2}{\partial t^2}
  I_{0,n}\left(\frac{\lambda}{v}
  \frac{\sqrt{n(n+1)}}{\sqrt[2n+2]{2n+2}}\sqrt[n+1]{y_1\ldots
  y_{n+1}}\right)dx_1\ldots dx_n.$$

  By analogy, for the $k$-th derivative we have
  ($k=\overline{3,n}$): $$\frac{\partial^k}{\partial t^k}\int\ldots\int_{T_{vt}}
  I_{0,n}\left(\frac{\lambda}{v}
  \frac{\sqrt{n(n+1)}}{\sqrt[2n+2]{2n+2}}\sqrt[n+1]{y_1\ldots
  y_{n+1}}\right)dx_1\ldots dx_n=$$ $$=\frac{(\sqrt{n+1})^{n+1}}{(\sqrt{n})^n(n-k)!}v^nt^{n-k}+
  \int\ldots\int_{T_{vt}}\frac{\partial^k}{\partial t^k}
  I_{0,n}\left(\frac{\lambda}{v}
  \frac{\sqrt{n(n+1)}}{\sqrt[2n+2]{2n+2}}\sqrt[n+1]{y_1\ldots
  y_{n+1}}\right)dx_1\ldots dx_n.$$

  Using the formulas found, we receive:

  $$\int\ldots\int_{T_{vt}}\left[\lambda^n+\lambda^{n-1}\frac{\partial}{\partial t}+\ldots+
  \frac{\partial^n}{\partial t^n}\right]I_{0,n}\left(\frac{\lambda}{v}
  \frac{\sqrt{n(n+1)}}{\sqrt[2n+2]{2n+2}}\sqrt[n+1]{y_1\ldots
  y_{n+1}}\right)dx_1\ldots dx_n=$$ $$=v^n\frac{(\sqrt{n+1})^{n+1}}{(\sqrt{n})^n}\left[e^{-\lambda t}-1-
  \lambda t-\ldots-\frac{(\lambda t)^{n-1}}{(n-1)!}\right].$$

  And finally, $$\frac{(\sqrt{n})^n}{(\sqrt{n+1})^{n+1}v^n}e^{-\lambda t}\int\ldots\int_{T_{vt}}\left[\lambda^n+\lambda^{n-1}\frac{\partial}{\partial t}+\ldots+
  \frac{\partial^n}{\partial
  t^n}\right]I_{0,n}\left(\frac{\lambda}{v}\times\right.
  $$ $$\left.\times
  \frac{\sqrt{n(n+1)}}{\sqrt[2n+2]{2n+2}}\sqrt[n+1]{y_1\ldots
  y_{n+1}}\right)dx_1\ldots dx_n = 1-e^{-\lambda t}-\ldots-\frac{(\lambda t)^{n-1}}{(n-1)!}e^{-\lambda t}.$$

  Theorem is proved.

  \underline{Remark 3.1} In the work [4] Orsingher and Sommella conjectured that the
  normalizing constant in (13) should be equal to $\frac{(\sqrt{n})^n}{(\sqrt{n+1})^{n+1}v^n}
  =\frac{t^n}{n!VolT_{vt}}.$ As we've seen it is not surprising
  because this constant appears when we find the integral $\int\ldots\int_{T_{vt}}
  1dx_1\ldots dx_n=Vol T_{vt}.$

\vskip 20mm

\centerline{References }

\vskip 2mm

[1] Samoilenko, I.V.(2001). Markovian random evolution in $R^n$.
\textit{Rand. Operat. and Stoc. Equat.}, N2, 139-160.

[2] Turbin, A.F. and Plotkin, D.Ja.(1991). Bessel equations of
high order. In: Asymptotic methods in problems of random
evolutions. \textit{Institute of Math. of NASU}, 112-121 (in
Russian).

[3] Orsingher, E.(2002). Bessel functions of third order and the
distribution of cyclic planar motions with three directions.
\textit{Stochastics and Stochastics Reports}, V74(3-4), 617-631.

[4] Orsingher, E. and Sommella, A.(2004). Cyclic random motion in
$R^3$ with four directions and finite velocity.
\textit{Stochastics and Stochastics Reports}, V76(2), 113-133.

\end{document}